\documentstyle[12pt]{article}
\input amssym.def
\input amssym
\begin{document}
\def \wt{{\rm wt}}
\def \spa{{\rm span}}
\def \End{{\rm End}}
\def \Hom{{\rm Hom}}
\def \Ind{{\rm Ind}}
\def \Mod{{\rm Mod}}
\def \Z{\Bbb Z}
\def \M{\Bbb M}
\def \C{\Bbb C}
\def \R{\Bbb R}
\def \Q{\Bbb Q}
\def \N{\Bbb N}
\def \ann{{\rm Ann}}
\def \min{{\rm min}}
\def \span{{\rm span}}
\def \<{\langle} 
\def \o{\omega}
\def \O{\Omega}
\def \M{{\cal M}}
\def \1t{\frac{1}{T}}
\def \>{\rangle} 
\def \t{\tau }
\def \a{\alpha }
\def \e{\epsilon }
\def \l{\lambda }
\def \L{\Lambda }
\def \g{\gamma}
\def \b{\beta }
\def \om{\omega }
\def \o{\omega }
\def \cg{\chi_g}
\def \ag{\alpha_g}
\def \ah{\alpha_h}
\def \ph{\psi_h}
\def \nor{\vartriangleleft}
\def \V{V^{\natural}}
\def \voa{vertex operator algebra\ }
\def \v{vertex operator algebra\ }
\def \be{\begin{equation}\label}
\def \bee{\begin{equation}\label}
\def \ee{\end{equation}}
\def \qed{\mbox{ $\square$}}
\def \pf {\noindent {\bf Proof:} \,}
\def \bl{\begin{lem}\label}
\def \el{\end{lem}}
\def \ba{\begin{array}}
\def \ea{\end{array}}
\def \bt{\begin{thm}\label}
\def \et{\end{thm}}
\def \ch{{\rm ch}}
\def \br{\begin{rem}\label}
\def \er{\end{rem}}
\def \ed{\end{de}}
\def \bp{\begin{prop}\label}
\def \ep{\end{prop}}
\def \Res{{\rm Res}}

\newtheorem{thm}{Theorem}[section]
\newtheorem{prop}[thm]{Proposition}
\newtheorem{coro}[thm]{Corollary}
\newtheorem{lem}[thm]{Lemma}
\newtheorem{rem}[thm]{Remark}
\renewcommand{\theequation}{\thesection.\arabic{equation}}

\begin{center}
{\large Twisted representations of vertex operator algebras and associative algebras}\\

Chongying Dong\footnote{Supported by NSF grant 
DMS-9303374 and a faculty research funds granted by the University of 
California at 
Santa Cruz.}\\
Department of Mathematics, University of
California, Santa Cruz, CA 95064 \\
Haisheng Li\footnote{Supported by NSF grant DMS-9616630.}\\
Department of Mathematical Sciences,
Rutgers University, Camden, NJ 08102\\
Geoffrey Mason\footnote{Supported by NSF grant DMS-9401272 and a faculty research funds granted by the University of 
California at 
Santa Cruz.}
\\
Department of Mathematics, University of
California, Santa Cruz, CA 95064
\end{center}
\hspace{1.5 cm}

\begin{abstract} Let $V$ be a vertex operator algebra and $g$ an 
automorphism of order $T$. We construct a sequence of associative
algebras $A_{g,n}(V)$ with $n\in\frac{1}{T}\Z$ nonnegative
such that $A_{g,n}(V)$ is a quotient of
$A_{g,n+1/T}(V)$ and a pair of functors between
the category of $A_{g,n}(V)$-modules which are not $A_{g,n-1/T}(V)$-modules
and the category of admissible $V$-modules. These
functors exhibit a bijection between the simple modules in each
category. We also show that $V$ is $g$-rational if and only if all 
$A_{g,n}(V)$ are finite-dimensional semisimple algebras.
\end{abstract}

\section{Introduction}
\setcounter{equation}{0}
In this paper we continue our study of twisted representations of vertex 
operator algebras and lay some further foundations of orbifold 
conformal field theory. 

Given a vertex operator algebra $V$ and an automorphism $g$ of finite
order $T,$ we have constructed an associative algebra $A_g(V)$ in [DLM1]
such that there is a bijection between simple $A_g(V)$-modules
and irreducible admissible $g$-twisted $V$-modules, generalizing 
Zhu's algebra $A(V)$ [Z]. It was proved in [DLM1] that $A_g(V)$ is
a finite-dimensional semisimple algebra if $V$ is $g$-rational. But
it is not clear whether the semisimplicity of $A_g(V)$ implies
the $g$-rationality of $V.$ 
Partially motivated by this in the case $g=1$  
and by the induced module theory for vertex operator algebra we 
constructed a series of associative algebras $A_n(V)$ for nonnegative integer
$n$ in [DLM2] so that
$A_{n-1}(V)$ is quotient algebra of $A_n(V)$ induced from the 
identity map on $V$ and that 
there is a one to one correspondence between the simple $A_n(V)$-modules
which cannot factor through $A_{n-1}(V)$ and irreducible admissible
$V$-modules. Moreover, $V$ is rational if and only if all $A_n(V)$ are 
finite-dimensional semisimple algebras. In the case $n=0$ we have 
Zhu's original algebra $A(V).$

In this paper we will study the twisted analogues of the algebras $A_n(V).$
In particular, we will construct a series of associative algebras
$A_{g,n}(V)$ for nonnegative numbers $n\in \frac{1}{T}\Z$ and show
that $A_{g,n-1/T}(V)$ is a natural quotient of $A_{g,n}(V).$ As in the
untwisted case, there is a bijection between the simple $A_{g,n}(V)$-modules
which cannot factor through $A_{g,n-1/T}(V)$ and the irreducible
admissible $g$-twisted $V$-modules. In the case $g=1$ we recover the algebras
$A_n(V)$ and the case $n=0$ amounts to the algebra $A_g(V).$  

Since most results in this paper are similar to those in [DLM2] where $g=1,$
we refer the reader in a lot of places of this paper to [DLM2] for
details. We assume that the reader is familiar with the elementary 
theory of vertex operator algebras as found in [B], [FLM], [FHL]
and the definition of twisted modules for vertex operator algebras and 
related ones as presented in [DLM1].  

This paper is organized as follows: In Section 2 we introduce
the associative algebras $A_{g,n}(V).$ In Section 3 we construct the functor
$\Omega_n$ from admissible $g$-twisted modules to $A_{g,n}(V)$-modules. 
We show that the each homogeneous subspaces of of the first $n$ pieces 
is a module for $A_{g,n}(V).$ The Section 4 is the heart of the paper.
In this section we construct the functor $L_n$ from $A_{g,n}(V)$-modules
to admissible $g$-twisted $V$-modules. The main strategy is to prove
the associativity for twisted vertex operators.

\section{The associative algebra $A_{g,n}(V)$}
\setcounter{equation}{0}
Let $V=(V,Y,{\bf 1},\omega)$ be a vertex operator algebra
and $g$ be an automorphism of $V$ of order $T.$ Then $V$ is a direct 
sum of eigenspaces of $T:$
\begin{equation}\label{2.1}
V=\bigoplus_{r=0}^{T-1}V^r
\end{equation}
where $V^r=\{v\in V|gv=e^{-2\pi ir/T}v\}.$ 

Fix $n=l+\frac{i}{T}\in\frac{1}{T}\Z$ with $l$ a nonnegative 
integer and $0\leq i\leq T-1.$ For $0\leq r\leq T-1$ we define 
$\delta_i(r)=1$ if $i\geq r$ and $\delta_i(r)=0$ if $i<r$.
We also set $\delta_i(T)=1.$ 
Let $O_{g,n}(V)$ be the linear span of all $u\circ_{g,n} v$ and $L(-1)u+L(0)u$
where for homogeneous $u\in V^r$ and $v\in V,$
\bee{g2.2}
u\circ_{g,n} v=\Res_{z}Y(u,z)v\frac{(1+z)^{\wt u-1+\delta_i(r)+l+r/T}}{z^{2l
+\delta_{i}(r)+\delta_{i}(T-r) }}.
\end{equation}

Define the linear space $A_{g,n}(V)$ to be the quotient $V/O_{g,n}(V).$ Then
$A_{g,n}(V)$ is the untwisted associative algebra $A_n(V)$ as defined in
[DLM2] if $g=1$ and is $A_g(V)$ in [DLM1] if $n=0.$

We also define a second product $*_{g,n}$ on $V$ for $u$ and $v$ as
above: 
\bee{a5.1}
u*_{g,n}v=\sum_{m=0}^{l}(-1)^m{m+l\choose l}\Res_zY(u,z)\frac{(1+z)^{\wt\,u+l}}{z^{l+m+1}}v
\end{equation} 
if $r=0$ and $u*_{g,n}v=0$ if $r>0.$
Extend linearly to obtain a bilinear product  on $V.$

\begin{lem}\label{l2.1} If $r\ne 0$ then $V^r\subset O_{g,n}(V).$
\end{lem}

\pf Let $v\in V^r$ be homogeneous. Then $v\circ_{g,n}{\bf 1}\in O_{g,n}(V).$
From the definition we know that
$$v\circ_{g,n}{\bf 1}=\sum_{j=0}^{\infty}{\wt v-1+l+\delta_{i}(r)+r/T\choose
j}v_{j-2l-\delta_i(r)-\delta_i(T-r)}{\bf 1}.$$
Note that $v_k{\bf 1}=0$ and $v_{-k-1}{\bf 1}=\frac{1}{k!}L(-1)^kv$ for $k\geq
0.$ Using $L(-1)v\equiv -L(0)v$ modulo $O_{g,n}(V)$ we see that
$v\circ_{g,n}{\bf 1}\equiv \sum_{j=0}^{2l-1+\delta_i(r)+\delta_i(T-r)}\frac{a_j}{j!}\frac{r^j}{T^j}v$
where $a_j$ are integers and $a_{2l-1+\delta_i(r)+\delta_i(T-r)}=\pm 1.$ 
Thus $v\circ_{g,n}{\bf 1}\equiv cv$ modulo $O_{g,n}$ for a nonzero constant
$c.$ This shows that $v\in O_{g,n}(V).$
\qed

\begin{lem}\label{l2.2} (i)
 Assume that $u\in V$ is homogeneous,
$v\in V$ and $m\ge k\ge 0.$ Then 
$$\Res_{z}Y(u,z)v\frac{(1+z)^{{\wt}u-1+l+\delta_i(r)+\frac{r}{T}+k}}{z^{2l+
\delta_i(r)+\delta_i(T-r)+m}}\in O_{g,n}(V).$$

(ii) For homogeneous $u,v\in V^0,$ 
$u*_nv-v*_nu-\Res_zY(u,z)v(1+z)^{\wt u-1}\in O_{g,n}(V).$
\end{lem}

\pf The proof of (i) is similar to that of Lemma 2.1.2 of [Z]. (ii) follows
from a result in Lemma 2.1 (iii) of [DLM2] that 
$u*_{g,n}v-v*_{g,n}u-\Res_zY(u,z)v(1+z)^{\wt u-1}\in O_{1,l}(V^0)$ and
the containment $O_{1,l}(V^0)\subset O_{g,n}(V).$ 
\qed

\begin{lem}\label{l2.3} (i) $O_{g,n}(V)$ is a 2 sided ideal of $V$ 
under $*_{g,n}.$

(ii) If $I=O_{g,n}\cap V^0$ then $I/O_{1,l}(V^0)$ is a two-sided ideal
of $A_{1,l}(V^0).$

\end{lem}

\pf Since $V^r$ ($r>0$) is a subset of $O_{g,n}$ by Lemma \ref{l2.1},
we see that $O_{g,n}(V)=I\oplus (\oplus_{r=1}^{T-1}V^r).$ Clearly
$V^0*_{g,n}V^r\subset V^r.$ So (i) and (ii) are equivalent. We prove (ii).
Choose $c\in V^0$ homogeneous and $u\in I.$ Using Lemma \ref{l2.2} (i)
and the argument used to prove Proposition 2.3 of [DLM1] we show that
both $c*_{g,n}u$ and $u*_{g,n}c$ lie in $O_{g,n}(V).$ 
\qed

The first main result is the following:
\begin{thm}\label{t2.4}  (i) The product $*_{g,n}$ induces the structure of an 
associative algebra  on $A_{g,n}(V)$ with identity ${\bf 1}+O_{g,n}(V).$

(ii) The linear map 
$$\phi:  v\mapsto e^{L(1)}(-1)^{L(0)}v$$
induces an anti-isomorphism $A_{g,n}(V)\to A_{g^{-1},n}(V)$.

(iii) $\omega+O_{g,n}(V)$ is a central element of $A_{g,n}(V).$
\end{thm}

\pf (i) follows from the result in [DLM2] that $A_{1,l}(V^0)$ is an
associative algebra with respect to $*_{g,n}$ and Lemma \ref{l2.3} (i).
The proof of (ii) is similar to that of Theorem 2.4 (ii) of [DLM1].
(iii) follows from Theorem 2.3 (iii) of [DLM2] which says that
 $\omega+O_{1,l}(V^0)$ is a central 
element of $A_{l,1}(V^0)$ and Lemma \ref{l2.3} (ii). \qed 

\begin{prop}\label{inver}
The identity map on $V$ induces an onto algebra homomorphism
from $A_{g,n}(V)$ to $A_{g,n-1/T}(V).$
\end{prop}

\pf If $n=l+i/T$ with $i\geq 1$ then both $A_{g,n}$ and $A_{g,n-1/T}$ are
quotients of $A_{1,l}(V^0).$ Otherwise $i=0$ and $A_{g,n}$
is a quotient algebra of $A_{1,l}(V^0)$ and $A_{g,n-1/T}$ is a quotient
algebra of $A_{1,l-1}(V^0).$ By Proposition 2.5 of [DLM2] the identity
map induces an epimorphism from $A_{1,l}(V^0)$ to $A_{1,l-1}(V^0).$ So
it is enough to show that $O_{g,n}(V)\cap V^0\subset O_{g,n-1/T}(V)\cap V^0.$
But this follows from Lemma \ref{l2.2} (i) immediately. \qed

As in [DLM2], Proposition \ref{inver} in fact gives us an 
inverse system $\{A_{g,n}(V)\}.$
Denote by $I_{g,n}(V)$ the inverse limit $\displaystyle{\lim_{\leftarrow}A_n(V).}$
Then\begin{equation}\label{i1}
I_g(V)=\{a=(a_{n}+O_{g,n}(V))\in \prod_{n\geq 0, n\in\frac{1}{T}\Z}A_{g,n}(V)| 
a_{n}-a_{n-1/T}\in O_{n-1/T}(v)\}.
\end{equation}
An interesting problem is to determine $I_g(V)$ explicitly and to
study the representations of $I_g(V).$

\section{The Functor $\Omega_n$}
\setcounter{equation}{0}
Recall from [DLM1] the Lie algebra $V[g]$
$$V[g]={\cal L}(V,g)/D{\cal L}(V,g)$$
where 
$${\cal L}(V,g)=\oplus_{r=0}^{T-1}t^{r/T}\C[t,t^{-1}]\otimes V^r$$
and by $D={d\over dt}\otimes 1+1\otimes L(-1).$ In order to write
down the Lie bracket we introduce the notation  $a(q)$ which
is  the image of $t^q\otimes a\in
{\cal L}(V,g)$ in $V[g].$ Let $a\in V^r,$ $v\in V^s$ and $m,n\in\Z.$ Then
$$[a(m+{r\over T}), b(n+{s\over T})]=
\sum_{i=0}^{\infty}{m+{r\over T}\choose i}a_ib(m+n+{r+s\over T}-i).$$
In fact $V[g]$
is $\frac{1}{T}\Z$-graded Lie algebra 
 by defining the degree of $a(m)$ to be $\wt v-m-1$ if $v$
is homogeneous. Denote the homogeneous subspace of degree $m$ by 
$\hat V[g]_m.$  In particular, $\hat V[g]_0$ is a Lie subalgebra. 

By Lemma \ref{l2.2} (ii) we have 
\bp{p2.5}
 Regarded $A_{g,n}(V)$ as a Lie algebra, the map $v(\wt v-1)\mapsto v+O_{g,n}(V)$ 
is a well-defined onto Lie algebra homomorphism from $\hat V[g]_0$ to
$A_{g,n}(V).$
\ep

For a module $W$ for the Lie algebra $V[g]$ and a nonnegative
$m\in \frac{1}{T}\Z$ we let $\O_{m}(W)$ denote the space
of ``m-th lowest
weight vectors,'' that is
\begin{eqnarray}\label{g5.1}
\Omega_m(W)=\{u\in W|V[g]_{-k}u=0\ if \ k\geq m\}.
\end{eqnarray}
Then $\O_{m}(W)$ is a module for the Lie algebra $V[g]_0.$ 

Note that 
for a $g$-twisted weak $V$-module $M$ the map $v(m)\mapsto v_m$ for $v\in V$ and $m\in \frac{1}{T}\Z$ gives a representation of $V[g]$ on $M$ [DLM2]. For
a homogeneous $v\in V$ we set $o_p(v)=v_{\wt v-1-p}$ on $M.$ 

\begin{lem}\label{kel} Let $M$ be a weak $V$-module.
Then for any homogeneous $u\in V^r,$ $v\in V^s,$  
$p\in \frac{r}{T}+Z, q\in \frac{s}{T}+Z$ with $p\geq q\geq -n$ and $p+q\geq 0$ 
these exists a unique $w_{u,v}^{p,q}\in V^{r+s}$ 
such that  $o_p(u)o_q(v)=o_{p+q}(w_{u,v}^{p,q})$ on  $\Omega_n(M).$ 
In particular if $s=T-r$ and $p=l+\delta_i(T-r)-k-\frac{r}{T}=-q$ 
for $k=0,...,l$  
$$w_{u,v}^{p,-p}=\sum_{m=0}^{k}(-1)^m{2l\!+\!\delta_i(r)\!+\!\delta_i(T-r)\!-\!1\!+\!m\!-\!k\choose m}\Res_zY(u,z)v\frac{(1+z)^{\wt u+l-1+\delta_i(r)+\frac{r}{T}}}{z^{2l+\delta_i(r)+\delta_i(T-r)-k+m}}.$$
\end{lem}  

The proof is similar to that of Theorem 3.2 of [DLM2]. This lemma is important
in constructing admissible $g$-twisted modules from $A_{g,n}(V)$-modules
in the next section.

\bt{t5.3} Suppose that $M$ is a weak $V$-module. Then there is 
a representation of the associative algebra $A_{g,n}(V)$ on $\O_{n}(M)$ induced by
the map $a\mapsto o(a)=a_{\wt a-1}$ for homogeneous $a\in V.$ 
\et

\pf By Theorem 5.1 of [DLM2],  the map $a\mapsto o(a)=a_{\wt a-1}$ for 
homogeneous $a\in V^0$ induces a representation of $A_{1,l}(V^0)$ on
$\O_{n}(M).$ So it is enough to show that $o(a)=0$ for
$a\in O_{g,n}(V)\cap V^0.$ It is clear that $o(L(-1)u+L(0)u)=0.$ 
It remains to show that $o(u\circ_{g,n}v)=0$ for $u\in V^r$ and $v\in V^{T-r}.$

Recall identity (10) from [DL]: for $ p \in \Z $ and $ s, t  \in \Q ,$  
\begin{equation}\label{jacobi-term}
\sum_{m\geq 0}(-1)^m {p\choose m}(u_{p+s-m}v_{t+m}-(-1)^{p}v_{p+t-m}u_{s+m})=
\sum_{m\geq 0}{s\choose m}(u_{p+m}v)_{s+t-m}.
\end{equation}
Now take $-p=2l+\delta_{i}(r)+\delta_{i}(T-r),$ $s=\wt u-1+l+\delta_i(r)
+\frac{r}{T},$ $t= \wt v-1+l+\delta_i(T-r)+\frac{T-r}{T}.$ Then on 
$\Omega_n(M)$ we have $u_{s+m}=v_{t+m}=0$ for $m\geq 0.$ Thus on $\Omega_n(M),$  
\begin{eqnarray*}
& & 0=\sum_{m\geq 0}{s\choose m}(u_{p+m}v)_{s+t-m}\\
& & \ \ =o(\Res_{z}Y(u,z)v\frac{(1+z)^{\wt u-1+\delta_i(r)+l+r/T}}{z^{2l
+\delta_{i}(r)+\delta_{i}(T-r) }})\\
& & \ \ =o(u\circ_{g,n}v).
\end{eqnarray*}
This completes the proof. \qed

Let $M=\oplus_{m\geq 0, m\in\frac{1}{T}\Z}M(m)$ be an admissible 
$g$-twisted module with $M(0)\ne 0.$ 

\bp{l2.9} The following hold

(i) $\Omega_n(M)\supset \oplus_{i=0}^{n}M(i).$ If $M$ is simple
then  $\Omega_n(M)=\oplus_{i=0}^{n}M(i).$

(ii) Each $M(p)$ is an $\hat V[g]_0$-module and
$M(p)$ and $M(q)$ are inequivalent if $p\ne q$ and both
$M(p)$ and $M(q)$ are nonzero.  If $M$ is simple then each
$M(p)$ is an irreducible $\hat V[g]_0$-module.

(iii) Assume that $M$ is simple. Then
each $M(i)$ for $i=0,...,n$ is a simple $A_{g,n}(V)$-module
and $M(i)$ and $M(j)$ are inequivalent $A_{g,n}(V)$-modules. 
\ep

The proof is similar to that of Proposition 3.4 of [DLM2].

\section{The functor $L_{n}$}
\setcounter{equation}{0}

In Section 3 we have shown how to obtain an $A_{g,n}(V)$-module from
an admissible $g$-twisted $V$-module. We show in this section that 
there is a universal
way to construct an admissible $g$-twisted $V$-module
from an $A_{g,n}(V)$-module which cannot factor through $A_{g,n-1/T}(V).$
(If it can factor through $A_{g,n-1/T}(V)$ we can consider the same
procedure for $A_{g,n-1/T}(V).$) As in [DLM2], 
a certain quotient of the universal object is an admissible $g$-twisted
$V$-module $L_n(U)$ and $L_n$ defines a functor which
is a right inverse to the functor $\O_n/\O_{n-1/T}$ 
where $\O_n/\O_{n-1/T}$ is the quotient functor $M\mapsto \O_n(M)/\O_{n-1/T}(M).$

Fix an $A_{g,n}(V)$-module $U$ which cannot factor through $A_{g,n-1/T}(V).$ Then
it is a  module for $A_{g,n}(V)_{Lie}$ in an obvious way.
By Proposition \ref{p2.5} we can lift $U$ to a module for the Lie algebra
$V[g]_0,$ and then to one for $P_n=\oplus_{p>n}V[g]_{-p}\oplus V[g]_0$ by letting $V[g]_{-p}$ act trivially. Define 
\bee{g6.1}
M_n(U)=\Ind_{P_n}^{V[g]}(U)=U(V[g])\otimes_{U(P_n)} U. 
\end{equation}
If we give $U$ degree $n$, the $\frac{1}{T}\Z$-gradation of $V[g]$ lifts to 
$M_n(U)$ which thus becomes a  $\frac{1}{T}\Z$-graded module for $V[g].$
The PBW theorem implies that $M_n(U)(i)=U(V[g])_{i-n}U.$ 

We define for $v\in V,$
\bee{g6.2}
Y_{M_n(U)}(v,z)=\sum_{m\in\frac{1}{T}\Z}v(m)z^{-m-1}
\end{equation}
As in [DLM1], $Y_{M(U)}(v,z)$ satisfies all conditions of a week $g$-twisted
$V$-module
except the associativity which does not hold on $M_n(U)$ in general.
We have to divide out by the desired relations. 

Let $W$ be the subspace of $M_n(U)$ spanned linearly by the 
coefficients of
\begin{eqnarray}\label{g6.3}
& &(z_{0}+z_{2})^{{\wt}a-1+l+\delta_i(r)+\frac{r}{T}}Y(a,z_{0}+z_{2})Y(b,z_{2})u\nonumber\\
& &\ \ \ \ \ \ -(z_{2}+z_{0})^{{\wt}a+-1+l+\delta_i(r)+\frac{r}{T}}
Y(Y(a,z_{0})b,z_{2})u
\end{eqnarray}
for any homogeneous $a\in V^r,b\in V,$ $u\in U$.
Set
\be{g6.4}
\bar M_n(U)=M_n(U)/U(V[g])W.
\end{equation}

\bt{t6.1} The space $\bar M_n(U)
=\sum_{m\geq 0}\bar M_n(U)(m)$ is an admissible $g$-twisted 
$V$-module with $\bar M_n(U)(0)\ne 0,$
$\bar M_n(U)(n)=U$ and with the 
following universal property: for any weak $g$-twisted $V$-module $M$
and any $A_{g,n}(V)$-morphism $\phi: U\to \O_n(M),$ there is a unique morphism 
$\bar\phi: \bar M_n(U)\to M$ of weak $g$-twisted $V$-modules which 
extends $\phi.$ 
\et

See Theorem 4.1 of [DLM2] for a similar proof.  

Let $U^*=\Hom_{\C}(U,\C)$ and let $U_s$ be the subspace
of $M_n(U)(n)$ spanned by  ``length'' $s$ vectors  
 $$o_{p_1}(a_1)\cdots o_{p_s}(a_s)U$$
where $p_1\geq \cdots \geq p_s,$ $p_1+\cdots p_s=0,$ $p_i\ne 0,$ $p_s\geq -n$
and $a_i\in V.$ Then by PBW theorem 
$M_n(U)(n)=\sum_{s\geq 0}U_s$ with $U_0=U$ and $U_s\cap U_t=0$ if $s\ne t.$
Motivated by the results in Lemma \ref{kel}  we extend $U^*$ to 
$M_n(U)(n)$ inductively so that 
\begin{eqnarray}\label{def}
\<u',o_{p_1}(a_1)\cdots o_{p_s}(a_s)u\>
=\<u',o_{p_{1}+p_2}(w_{a_1,a_2}^{p_1,p_2})o_{p_3}(a_3)\cdots o_{p_{s}}(a_{s})u).
\end{eqnarray}
where $o_j(a)=a(\wt a-1-j)$ for homogeneous $a\in V.$ 
We further extend
$U^*$ to $M_n(U)$ by letting $U^*$ annihilate $\oplus_{i\ne n}M(U)(i).$

Set 
$$ J=\{v\in M_n(U)|\langle u',xv\rangle=0\ {\rm for\ all}\ u'\in
U^{*},\ {\rm all}\ x\in U(V[g])\}.$$
We can now state the second main result of this section.
\bt{t6.3} Space $L_n(U)\!=\!M_n(U)/J$ is an admissible
$V$-module satisfying $L_n(U)(0)\ne 0$ and
$\O_n/\O_{n-1/T}(L_n(U))\cong U.$
Moreover 
$L_n$ defines a functor from the category of  
$A_{g,n}(V)$-modules which cannot factor through $A_{g,n-1/T}(V)$  to
the category of admissible $V$-modules such that $\O_n/\O_{n-1/T}\circ L_n$
is naturally equivalent to the identity. Moreover, $L_n(U)$ is
a quotient module of $\bar M_n(U).$ 
\et

The proof of this theorem is the most complicated one. Fortunately we
can following the proof of Theorem 4.2 of [DLM2] step by step with suitable
modifications.  Again we refer the reader to [DLM2] for details. 

The analogue of Theorem 4.9 of [DLM2] (whose proof is easy) 
is the following.

\bt{t7.2} $L_n$ and $\O_n/O_{n-1/T}$
are equivalences when restricted to the full 
subcategories of completely reducible $A_{g,n}(V)$-modules whose irreducible
components cannot factor through $A_{g,n-1/T}(V)$ 
and 
completely reducible admissible $g$-twisted $V$-modules respectively.
In particular, $L_n$ and $\Omega_n/\O_{n-1/T}$ induce mutually
inverse bijections on the isomorphism classes of simple objects
in the category of $A_{g,n}(V)$-modules which cannot factor through $A_{g,n-1/T}(V)$
and admissible $g$-twisted $V$-modules
respectively.
\end{thm}

We also have the generalization of Theorem 4.10 of [DLM2] with a similar
proof. 

\bt{t8.1} Suppose that $V$ is a $g$-rational vertex operator algebra. Then the 
following hold:

(a) $A_{g,n}(V)$ is a finite-dimensional, semisimple associative algebra.

(b) The functors $L_n$ and $\O_n/O_{n-1/T}$ are mutually inverse categorical equivalences
between the category of $A_{g,n}(V)$-modules whose irreducible components
cannot factor through $A_{g,n-1/T}(V)$ and the category of admissible 
$g$-twisted $V$-modules.

(c) The functors $L_n,\O_n/\O_{n-1/T}$ induce mutually inverse categorical equivalences
between the category of finite-dimensional 
$A_{g,n}(V)$-modules whose irreducible components
cannot factor through $A_{g,n-1/T}(V)$ and the category of ordinary $g$-twisted
$V$-modules.
\et

As in [DLM2] one expects that if $A_{g,0}(V)$ is semisimple 
then $V$ is $g$-rational. We present some
partial results which are applications of $A_{g,n}(V)$-theory. 
\begin{thm}\label{t4.13} All $A_{g,n}(V)$ are finite-dimensional 
semisimple algebras if and only if $V$ is $g$-rational.
\end{thm}


\begin{thebibliography}{FLM1}


\bibitem[B]{B}
R. E. Borcherds, Vertex algebras, Kac-Moody algebras, and the Monster,
{\it Proc. Natl. Acad. Sci. USA} {\bf 83} (1986), 3068-3071.

\bibitem[DLM1]{DLM1} C. Dong, H. Li and G. Mason, Twisted representations of 
vertex operator algebras, q-alg/9509005.

\bibitem[DLM2]{DLM2} C. Dong, H. Li and G. Mason,
Vertex operator algebras and associative algebras, q-alg/9612010.

\bibitem[DL]{DL} C. Dong and Z. Lin, Induced modules for vertex operator algebras, {\em Commu. Math. Phys.}
{\bf 179} (1996), 157-184.

\bibitem[FHL]{FHL}
I. Frenkel, Y. Huang and J. Lepowsky, On axiomatic approaches to
vertex operator algebras and modules, Memoirs Amer. Math.
Soc. {\bf 104}, 1993.

\bibitem[FLM]{FLM}
I. Frenkel, J. Lepowsky and A. Meurman, {\it Vertex Operator Algebras
and the Monster}, Pure and Appl. Math., {\bf Vol. 134}, Academic Press,
Boston, 1988.

\bibitem[Z]{Z}
Y. Zhu, Modular invariance of characters of vertex operator algebras,
{\em J. Amer, Math. Soc.} {\bf 9} (1996), 237-302.
\end{thebibliography}
\end{document}